\documentclass[12pt]{amsart}
\usepackage{amssymb,amsfonts,amsmath,amsopn,amstext,amscd,latexsym, amsthm, enumerate,mathrsfs}

\usepackage{color}
\usepackage{xcolor}
\usepackage{hyperref}
\hypersetup{
     colorlinks   = true,
     citecolor    = blue
}

\usepackage[margin=30mm]{geometry}
\usepackage{bbm}
\usepackage[all,cmtip]{xy}

\usepackage{enumerate}
\usepackage[shortlabels]{enumitem}

\newtheorem{theorem}{Theorem}[section]

\newtheorem{lemma}[theorem]{Lemma}

\newtheorem{corollary}[theorem]{Corollary}

\theoremstyle{remark}

\DeclareMathOperator{\Aut}{Aut}

\numberwithin{equation}{section}

\newcommand{\Cay}{{\mathrm {Cay}}}

\newcommand{\SC}{{\mathrm {SC}}}
\newcommand{\eqd}{$\hfill \square$\medskip}

\begin{document}

\title[Normality of one-matching semi-Cayley graphs ]{Normality of one-matching semi-Cayley graphs over finite abelian groups with maximum degree three}

\author[M.~Arezoomand]{Majid Arezoomand}
\address{University of Larestan, Larestan, 74317-16137, Iran}

\email{arezoomand@lar.ac.ir(Corresponding author)}

\author[M.~Ghasemi]{Mohsen Ghasemi}
\address{Department of Mathematics, Urmia University,  Urmia  57135, Iran}
\email{m.ghasemi@urmia.ac.ir}

\thanks{}





\begin{abstract}
A graph $\Gamma$ is said to be a semi-Cayley graph over a group $G$
if it admits $G$ as a semiregular automorphism group with two orbits
of equal size. We say that $\Gamma$ is normal if $G$ is a normal
subgroup of ${\rm Aut}(\Gamma)$. We prove that 
every connected intransitive one-matching semi-Cayley graph, with maximum degree three,
 over a finite abelian group is normal and characterize all such non-normal graphs.
\end{abstract}

\maketitle
\section{Introduction}

Throughout this paper, groups are finite, and graphs
are finite, connected, simple and undirected. For the
graph-theoretic and group-theoretic terminology not defined here, we
refer the reader to ~\cite{Biggs, Wi}. Let $G$ be a
permutation group on $\Omega$ and $\alpha \in \Omega$. Denote by
$G_\alpha$ the stabilizer of $\alpha$ in $G$, that is, the subgroup
of $G$ fixing the point $\alpha$. We say that $G$ is semiregular on
$\Omega$ if $G_\alpha=1$ for every $\alpha \in \Omega$ and regular
if $G$  is transitive and semiregular.
Let $G$ be a group and $S$ a subset of $G$ not containing the
identity element $1_G$. The Cayley digraph $\Gamma= {\rm Cay(G, S)}$
of $G$ with respect to $S$ has vertex set $G$ and arc set $\{(g,
sg)\mid g \in G, s \in S\}$. If $S=S^{-1}$ then ${\rm Cay(G, S)}$
can be viewed as an undirected graph, identifying an undirected edge
with two directed edges $(g, h)$ and $(h, g)$. This graph is called
Cayley graph of $G$ with respect to $S$. By a theorem of
Sabidussi ~\cite{S}, a graph $\Gamma$ is a Cayley graph over a group $G$
if and only if there exists a regular subgroup of ${\rm
Aut(\Gamma)}$ isomorphic to $G$.

There is a natural generalization of the Sabidussi's Theorem. A
graph $\Gamma$ is called an $n$-Cayley graph over a group $G$ if
there exists an $n$-orbit semiregular subgroup of ${\rm
Aut(\Gamma)}$ isomorphic to $G$. Undirected and loop-free $2$-Cayley graphs are called {\it semi-Cayley} \cite{AT,resmini}, and also
{\it bi-Cayley} by some authors \cite{Zhou}. 
$n$-Cayley graphs have been playing an
important role in many classical fields of graph theory, such as
strongly regular graphs  ~\cite{Le, M, Mar, resmini}, Hamiltonian graphs ~\cite{W}
$n$-extendable graphs ~\cite{G1,  Lu}, the spectrum of graphs ~\cite{arezoomand1, arezoomand3, G2},  automorphisms
\cite{arezoomand, arezoomand2, HKM, Zhou}, and the
connectivity of graphs ~\cite{C, Li}.

A graph $\Gamma$ is called a {\it semi-Cayley graph} over a group $G$ if $\Aut(\Gamma)$ admits a semiregular subgroup
 $R_G$ isomorphic to $G$ with two orbits (of equal size). Let $\Gamma$ be a semi-Cayley graph over a group $G$. Then there exists subsets
 $R,L$ and $S$ of $G$ such that $R=R^{-1}$, $L=L^{-1}$ and $1\notin R\cup L$ such that $\Gamma\cong\SC(G;R,L,S)$, where
 $\SC(G;R,L,S)$ is an undirected graph with vertices $G\times\{1,2\}$ and its edge set consists of three sets (see \cite[Lemma 2.1]{resmini}):
 \begin{eqnarray*}
&&\{\{(x,1),(y,1)\}\mid yx^{-1}\in R\}\ \ \  (\textrm{right edges}),\\
&&\{\{(x,2),(y,2)\}\mid yx^{-1}\in L\}\ \ \   (\textrm{left edges}),\\
&&\{\{(x,1),(y,2)\}\mid yx^{-1}\in S\} \ \ \  (\textrm{spoke edges}).
\end{eqnarray*}
Furthermore, $R_G:=\{\rho_g\mid g\in G\}$, where
$\rho_g:G\times\{1,2\}\rightarrow G\times\{1,2\}$ and
$(x,i)^{\rho_g}=(xg,i)$, $i=1,2$, is a semiregular subgroup of
$\Aut(\SC(G;R,L,S))$ isomorphic to $G$ with two orbits
$G\times\{1\}$ and $G\times\{2\}$.
 A semi-Cayley graph
$\Gamma=\SC(G;R,L,S)$ over a group $G$ is called {\it normal}
over $G$ if $R_G$ is a normal subgroup of $\Aut(\Gamma)$ (see \cite[p. 42]{arezoomand}) and it
 is called {\it one-matching} if $S=\{1\}$ (see \cite[p.
603]{Kovacs}).
In this paper, we prove:

\begin{theorem}\label{asli}
Let $\Gamma=\SC(G;R,L,\{1\})$ be a connected one-matching semi-Cayley graph over a finite abelian group $G\neq 1$ with $|R|,|L|\leq 2$.
Then $\Gamma$ is normal if and only if  none of the following are satisfied (even after interchanging $R$ and $L$)
\begin{itemize}
\item[$(1)$] $|R|=|L|=1$ (so $G\cong\Bbb Z_2$ or $\Bbb Z_2^2$),
\item[$(2)$] $|R|=|L|=2$ and $|R\cap L|=1$ (so $G\cong\Bbb Z_2^2$ or $\Bbb Z_2^3$),
\item[$(3)$] $R=L=\{a,a^{-1}\}$, where $o(a)=4$ (so $G=\langle a\rangle\cong\Bbb Z_4$),
\item[$(4)$] $R=\{a,b\}$, $L=\{c,c^{-1}\}$, where $o(a)=o(b)=2$, $o(c)=4$ and
$G=\langle a\rangle\times\langle b\rangle\times\langle c\rangle\cong\Bbb Z_2^2\times\Bbb Z_4$,
\item[$(5)$] $R=\{a,a^{-1}\}$, $L=\{b,ba^2\}$, where $o(a)=4$, $o(b)=2$ and $G=\langle a\rangle\times\langle b\rangle\cong\Bbb Z_4\times\Bbb Z_2$,
\item[$(6)$] $R=\{a,a^{-1}\}$, $L=\{a^k,a^{-k}\}$, where $o(a)=n$ and $(n,k)$ is one of the pairs $(5,2)$, $(8,3)$, $(10,2)$,
$(10,3)$, $(12,5)$
or $(24,5)$
(so $G\cong\Bbb Z_n$),
\item[$(7)$] $R=\{a,a^{-1}\}$, $L=\{a^3b,a^{-3}b\}$ or $L=\{a^2b,a^{-2}b\}$, where $o(a)=10$, $o(b)=2$ and $G=\langle a\rangle\times\langle b\rangle \cong\Bbb Z_{10}\times\Bbb Z_2$,
\item[$(8)$] $R=\{a,a^{-1}\}$, $L=\{ab, a^{-1}b\}$, where $o(a)=4$, $o(b)=2$ and $G=\langle a\rangle\times\langle b\rangle \cong\Bbb Z_4\times\Bbb Z_2$.
\end{itemize}
Furthermore, in all of the above cases, $\Gamma$ is transitive.
\end{theorem}

For a graph $\Gamma $, we use  $V(\Gamma )$, $E(\Gamma )$,  $A(\Gamma )$ and  ${\rm Aut(\Gamma )}$
to denote its vertex set, edge set,  arc set and its  full automorphism group
respectively. For $v \in V (\Gamma )$,  $N(u)$ is the neighborhood of $u$ in $\Gamma $, that is, the set of vertices adjacent
to $u$ in $\Gamma $.  A graph $\Gamma$ is called transitive if $\Aut(\Gamma)$ is
transitive on $V(\Gamma)$, otherwise it is called intransitive.
Also  a graph $\Gamma $ is said to be  edge-transitive and arc-transitive (or
symmetric) if Aut$(\Gamma )$ acts transitively on  E$(\Gamma )$  and A$(\Gamma )$, respectively.

\section{Preliminaries}

Let $\Gamma=\SC(G;R,L,\{1\})$ be a one-matching semi-Cayley graph over a finite group $G\neq 1$.
Let $\Gamma_0=\SC(G;L,R,\{1\})$ the graph obtained from interchanging the left and right edges of $\Gamma$. Then
 $\Gamma\cong\Gamma_0$. Furthermore, $\Aut(\Gamma)\cong\Aut(\Gamma_0)$ and also $R_G\unlhd\Aut(\Gamma)$
 if and only if $R_G\unlhd\Aut(\Gamma_0)$. Hence, in studying the normality of $\Gamma$, we may assume that
 $|L|\leq |R|$. Moreover, since $\Gamma$ is a normal over a group $G$ if and if its
 complement $\Gamma^c$ is normal over $G$, we may assume that $\Gamma$ is connected or equivalently $G=\langle R\cup L\rangle$.


Let $\Gamma=\SC(G;R,L,\{1\})$ be a connected semi-Cayley graph over a finite abelian group $G$, and let $A$ and $V$, be its automorphism group
and vertex set, respectively. For each $\sigma\in\Aut(G)$ we
define two maps
\begin{eqnarray*}
\varphi_\sigma&:&V(\Gamma)\rightarrow V(\Gamma);~~~ (x,1)^{\varphi_\sigma}=(x^\sigma,1),~(x,2)^{\varphi_\sigma}=(x^\sigma,2),\\
\psi_\sigma&:&V(\Gamma)\rightarrow V(\Gamma);~~~(x,1)^{\psi_\sigma}=(x^\sigma,2),~(x,2)^{\psi_\sigma}=(x^\sigma,1).
\end{eqnarray*}
Set
\[X:=\{\varphi_\sigma\mid \sigma\in\Aut(G), R^\sigma=R, L^\sigma=L\},\]
\[Y:=\{\psi_\sigma\mid \sigma\in\Aut(G), R^\sigma=L, L^\sigma=R\},\]
 and let us denote $X\cup Y$ by $\Aut(G;R,L)$.
Then $N_A(R_G)=R_G\rtimes \Aut(G;R,L)$ by \cite[Theorem 1]{arezoomand}. So
$R_G\unlhd A$ if and only if $A=R_G\rtimes \Aut(G;R,L)$ \cite[Proposition 2 (1)]{arezoomand}. Moreover, if $R_G\unlhd A$, then $A_{(1,1)}=X$
and the converse holds if $\Gamma$ is intransitive \cite[Proposition 2 (2)]{arezoomand}.
Also if $R_G\unlhd A$ then $\Gamma$ is intransitive if and only if $A_{(1,1)}=\Aut(G;R,L)$ \cite[Corollary 2.9]{arezoomand}.
Note that if $Y\neq\emptyset$, then $\Gamma$ is transitive. So if $R_G\unlhd A$ then
$\Gamma$ is transitive if and only if $Y\neq\emptyset$. Also, by the following lemma and above results, if $\Gamma$ is intransitive or $Y\neq\emptyset$, then $\Gamma$ is normal if and only if $A_{(1,1)}=X$. It is easy to see that
$A_{(1,1)}\cap N_A(R_G)=A_{(1,2)}\cap N_A(R_G)=X$. In particular, if $R_G\unlhd A$ then $A_{(1,1)}=A_{(1,2)}=X$.
In what follows, unless otherwise stated, we keep the above notations and use the above results without referring them.

\begin{lemma}\label{non-empty} Let $Y\neq\emptyset$.
Then $\Gamma$ is normal if and only if $A_{(1,1)}=X$.
\end{lemma}
\proof If $\Gamma$ is normal then $A_{(1,1)}=X$. Conversely, suppose that $A_{(1,1)}=X$.  Let $\beta\in A$ be arbitrary. We have to show
that $\beta\in N_A(R_G)$. Since $Y\neq\emptyset$ (and $Y\subseteq N_A(R_G)$), we may assume that $(1,1)^\beta\in G\times\{1\}$ (if $(1,1)^\beta\in G\times\{2\}$, then we replace $\beta$ with $\beta y$ for some $y\in Y$). Then after multiplying by an element of
$R_G$, we may assume that $(1,1)^\beta=(1,1)$. So $\beta\in A_{(1,1)}=X\subseteq N_A(R_G).$ \eqd

\section{Proof of Theorem \ref{asli}}
Keeping the notations of previous section, recall that $\Gamma=\SC(G;R,L,\{1\})$ is a connected semi-Cayley graph
over a finite abelian group $G\neq 1$ with $|L|\leq |R|\leq 2$, and $A$ denotes the automorphism group of $\Gamma$.
To prove Theorem \ref{asli}, we consider the all possibilities for the orders of $R$ and $L$ and their intersection.

Let us start with the following lemma:

\begin{lemma}\label{arc} Let $\Gamma$ be edge-transitive. Then it is non-normal. Also if $\Gamma$ is arc-transitive then $\Gamma$
is non-normal.
\end{lemma}
\proof It is enough to note that any element of the normalizer of $R_G$ must map $G$-orbits to $G$-orbits but an element
of $A$ that takes a right edge or left edge to a spoke edge does not do this. Since every connected arc-transitive graph is edge-transitive,
the second part is clear.
\eqd


\begin{lemma}\label{3.3} Let $L=\emptyset$, $R\neq\emptyset$. Then $\Gamma$ is intransitive and normal, and\\
$(1)$ if $|R|=1$ then $G\cong\mathbb Z_2$, $\Gamma\cong P_4$, $A\cong\Bbb Z_2$,\\
$(2)$ if $|R|=2$  then $G\cong\Bbb Z_n$ or $\Bbb Z_2\times\Bbb Z_2$, $n\geq 3$, and $A\cong D_{2|G|}$.
\end{lemma}
\proof $(1)$ It is clear.

$(2)$ Since $\Gamma$ is connected and $L=\emptyset$,  we have $G=\langle R \rangle\cong\Bbb Z_n$ or $\Bbb Z_2\times\Bbb Z_2$, for some $n\geq 3$. Hence $\Cay(G,R)$ is a $|G|$-cycle. By
\cite[Lemma 4.1]{arezoomand}, $A\cong\Aut(\Cay(G,R))\cong D_{2|G|}$.
\eqd

\begin{lemma}\label{3.4} Let $R=\{a\}$ and $L=\{b\}$. Then $\Gamma$ is transitive and non-normal and one of the following holds:\\
$(1)$ $G\cong\Bbb Z_2$, $A\cong D_8$.\\
$(2)$ $G\cong\Bbb Z_2\times\Bbb Z_2$, $A\cong D_{16}$.
\end{lemma}
\proof If $a=b$ then $G\cong\Bbb Z_2$ and otherwise $G\cong\Bbb Z_2\times\Bbb Z_2$. In both cases, $\Gamma$ is a $2|G|$-cycle and
so $A\cong D_{4|G|}$. Furthermore, in both cases $A_{(1,1)}\neq X$, which implies that both are non-normal.
\eqd

\begin{lemma}\label{zhou} Let $\Gamma$ be intransitive, $R\cap L=\emptyset$ and $\Gamma_\Omega$ be the quotient graph of $\Gamma$ with respect
to the one-matching set $\Omega=\{\{(g,1),(g,2)\}\mid g\in G\}$.
Then $A\leq\Aut(\Gamma_\Omega)$, where $\Gamma_\Omega$ is a
Cayley graph of $R_G$ with respect to $S=\{\rho_r,\rho_l\mid r\in R,
l\in L\}$ of valency $|R|+|L|$. In particular, if $\Gamma_\Omega$ is
a normal Cayley graph of $R_G$ then $\Gamma$ is a normal semi-Cayley
graph of $R_G$.
\end{lemma}
\proof We consider the action of $A$ on $\Omega$.
Let $K$ be the kernel of this action. Since $\Gamma$ is
intransitive, it implies that $K=1$ and so $A\leq {\rm
Aut(\Gamma_\Omega)}$. Clearly $R_G$ acts transitively on ${\rm
V(\Gamma_\Omega)}$. Now suppose that $\rho_h \in R_G$
 and $\{(g, 1), (g, 2)\}^{\rho_h}=\{(g, 1), (g, 2)\}$. Therefore $(g, 1)^{\rho_h}=(g, 1)$ and $(g, 2)^{\rho_h}=(g, 2)$ and so
 $(gh, 1)=(g, 1)$. Thus $\rho_h=1$ and $R_G$ acts regularly on ${\rm
V(\Gamma_\Omega)}$ and so $\Gamma_\Omega$ is a Cayley graph on $R_G$  with respect to $S$. Also since $R\cap L=\emptyset$, it
implies that $\Gamma_\Omega$ has valency $|R|+|L|$. \eqd

\begin{lemma}\label{3.6} If $|R|=2$  and $|L|=1$ then $\Gamma$ is normal.
\end{lemma}
\proof Let $R=\{a,b\}$ and $L=\{c\}$. If $c=a$ or $c=b$ then $a^2=b^2=1$, $A=R_G\cong\Bbb Z_2^2$ and so
$\Gamma$ is normal. Hence, we may assume that $c\neq a,b$. Suppose, towards a contradiction, that $\Gamma$ is non-normal.
Then $R\cap L=\emptyset$. Let $\Omega=\{\{(g,1),(g,2)\}\mid g\in G\}$
and $\Gamma_\Omega$ be the Cayley graph of $R_G$ with respect to $S=\{\rho_a,\rho_b,\rho_c\}$. Since
$\Gamma$ is non-normal, Lemma \ref{zhou} and \cite[Theorem 1.2]{Baik} imply that one of the following happens:\\
$(i)$ $o(a)=4$, $b=a^{-1}$ and $c=a^2$.\\
$(ii)$ $o(a)=4$, $b=a^{-1}$, $c^2=1$ and $c\notin\langle a\rangle$.\\
$(iii)$ $o(a)=6$, $b=a^{-1}$ and $c=a^3$.\\
In the first case, $A\cong D_8$ and $\Gamma$ is normal, in the second case $A\cong\Bbb Z_2\times D_8$ and $\Gamma$ is normal,
and in the last case, $A\cong D_{12}$ and $\Gamma$ is normal. Hence we get a contradiction.\eqd

\begin{lemma}\label{3.7} Let $R=L$, $|R|=2$. Then $\Gamma$ is transitive and the following are equivalent:\\
$(1)$ $\Gamma$ is normal.\\
$(2)$ $\Gamma$ is not arc-transitive.\\
$(3)$ $R=\{a,a^{-1}\}$, where $a$ is of order $k>2$ and $k\neq 4$.
\end{lemma}
\proof It is easy to see that $\Gamma$ is isomorphic to the $n$-prism graph, the cartesian product of an $n$-cycle with a path with
two vertices, where $n=|G|$, which is isomorphic to a Cayley graph on the dihedral group $D_{2n}=\langle s,t\mid s^n=t^2=(st)^2=1\rangle$, with respect to $S=\{s,s^{-1},t\}$. Hence $\Gamma$ is transitive.

By Lemma \ref{arc}, (1) implies (2). Now suppose that (2)
holds. If $R=\{b,c\}$, where $b^2=c^2=1$, then $\Gamma$ is
isomorphic to the three dimensional hypercube, which is
arc-transitive, a contradiction. Hence $R=\{a,a^{-1}\}$, where $a$
is of order $k>2$. Hence $G=\langle a\rangle\cong\Bbb Z_k$. Hence,
by \cite[Theorem 1.1]{Kovacs}, $k\neq 4$. Thus (2) implies (3). To
complete the proof, it is enough to prove that (3) implies (1). Let
(3) holds. Then $G\cong\Bbb Z_k$ and it is easy to see that $\Gamma$
is isomorphic to generalized Petersen graph $ GP(k,1)$ (see ~\cite{F}).
Also  by  \cite[Theorems 1 and
2]{F}  $GP(k,1)$ is vertex transitive and  $A\cong D_{2k}\times\Bbb Z_2$.
Hence $\Gamma$ is
vertex transitive and so $|A|=|A_{(1,1)}|2k$. This shows that
$|A_{(1,1)}|=2$. Since $R=L=\{a,a^{-1}\}$, $Y\neq\emptyset$ and
$|X|\geq 2$. Since $X\leq A_{(1,1)}$, we have $X=A_{(1,1)}$. Hence $\Gamma$ is normal, i.e. (1) holds.
This completes the proof. \eqd

\begin{lemma}\label{3.8} Let $|R|=|L|=2$, $|R\cap L|=1$. Then $\Gamma$ is transitive and non-normal. Also one of the following holds:\\
$(1)$ $G=\langle a,b\rangle\cong\Bbb Z_2\times\Bbb Z_2$, $R=\{a,b\}$ and $L=\{ab,b\}$.\\
$(2)$ $G=\langle a,b,c\rangle\cong\Bbb Z_2\times\Bbb Z_2\times\Bbb Z_2$, $R=\{a,b\}$ and $L=\{b,c\}$.
\end{lemma}
\proof
Since $R=R^{-1}$ and $L=L^{-1}$, both $R$ and $L$ consist of two involutions. Assume that $R=\{a,b\}$ and $L=\{b,c\}$.
 Since $G=\langle a,b,c\rangle$, if $c=ab$, $G=\langle a,b\rangle\cong\Bbb Z_2\times\Bbb Z_2$, otherwise $G\cong\Bbb Z_2\times\Bbb Z_2\times\Bbb Z_2$.
  In the former case,  \[\sigma=((ab, 1), (b, 2))((a, 2), (ab, 2))((a, 1), (1, 2)) \in A_{(1, 1)}\] but $\sigma \notin X$. Therefore
  $\Gamma$ is not normal.
    In the latter,
    \[\sigma=((a, 1), (1, 2))((a, 2), (c, 2))
    ((ab, 1), (b, 2))((ab, 2), (bc, 2))((ac, 2), (c, 1)) \in A_{(1, 1)}\]
     but $\sigma \notin X$. So $\Gamma$ is not normal. Also in both cases we see that
    $\Gamma$ is transitive.
\eqd

\begin{lemma}\label{2-4}
Let $|R|=|L|=2$, $R\cap L=\emptyset$. If $R=\{a,b\}$, where $a^2=b^2=1$, then one of the following holds:\\
$(1)$ $L=\{ab,c\}$, where $c^2=1$. In this case $G\cong\Bbb Z_2^3$, $\Gamma$ is intransitive and normal.\\
$(2)$ $L=\{c,d\}$, where $c^2=d^2=1$. In this case $G\cong\Bbb Z_2^4$ and $\Gamma$ is transitive and normal.\\
$(3)$ $L=\{c,c^{-1}\}$, where $c$ is of order $n>2$. In this case,
$G\cong\Bbb Z_2^2\times\Bbb Z_n$, and $\Gamma$
is normal if and only if $\Gamma$ is intransitive if and only if $n\neq 4$.\\
$(4)$ $L=\{c,c^{-1}\}$, where $o(c)=n>2$ is even, and $b=c^{n/2}$. In this case, $G\cong\Bbb Z_n\times\Bbb Z_2$, and $\Gamma$ is
normal and intransitive.\\
$(5)$ $L=\{c,c^{-1}\}$, where $o(c)=n>2$ is even, and $b=ac^{n/2}$. In this case, $G\cong\Bbb Z_n\times\Bbb Z_2$ and
$\Gamma$ is normal if and only if $\Gamma$ is intransitive if and only if $n\neq 4$.
\end{lemma}
\proof It is obvious that the possibilities of $L$ are exactly the
same given in (1)-(5).

$(1)$ In this case,
$G=\langle a,b,c\rangle\cong\Bbb Z_2^3$. Then, by GAP \cite{Gap}, $\Gamma$ is intransitive and normal.

 $(2)$ In this case, $G\cong\Bbb Z_2^4$
and by GAP , $A\cong (D_8\times D_8)\rtimes\Bbb Z_2$, $\Gamma$ is
transitive and normal.

$(3)$ Suppose that $L=\{c,c^{-1}\}$, where $c$ is an element of order $n>2$. Then $G\cong\Bbb Z_2^2\times\Bbb Z_n$. We prove that $\Gamma$ is normal if and only if it is intransitive if and only if $n\neq 4$.

If $n=4$ then by GAP, $\Gamma$ is transitive.  Conversely, suppose that $\Gamma$ is transitive. Then there exists
$\alpha\in\Aut(\Gamma)$ such that $(1,1)^\alpha=(1,2)$. Then $\alpha$ maps the $4$-cycle
\[(1,1),(b,1),(ab,1),(a,1),(1,1)\]
to a $4$-cycle including the point $(1,2)$. Since $R\cap L=\emptyset$, we have $(a,1)^\alpha,(b,1)^\alpha\neq (1,1)$. Hence
$(a,1)^\alpha,(b,1)^\alpha\in\{(c,2),(c^{-1},2)\}$, which implies that $(ab,1)^\alpha=(c^2,2)=(c^{-2},2)$. This means that $n=4$.

Let $\Gamma$ is normal. Then, since $Y=\emptyset$,  $\Gamma$ is intransitive. Conversely, suppose that $\Gamma$ is intransitive. So
$n\neq 4$, by the above discussion. Now \cite[Theorem 1.2]{Baik} and Lemma \ref{zhou}, imply that $\Gamma$ is normal.

$(4)$
In this case $G\cong\Bbb Z_n\times\Bbb Z_2$.
If $n=4$, then by GAP, $\Gamma$ is intransitive and normal. Hence, we may assume that $n\neq 4$.
So, by a similar argument of the previous case, $\Gamma$ is intransitive.
Suppose, towards a contradiction, that $\Gamma$ is non-normal. Then, by Lemma \ref{zhou} and \cite[Theorem 1.2]{Baik},
$n=6$. Now, by GAP, $A\cong\Bbb Z_2^2\times S_3$ which implies that $\Gamma$ is normal, a contradiction.

$(5)$ In this case $G\cong\Bbb Z_n\times\Bbb Z_2$.
If $n=4$ then $\Gamma$ is transitive and non-normal. Let $n\neq 4$. Then by
a similar argument of the case $(3)$, $\Gamma$ is intransitive.  By the same argument in case $(4)$, if $\Gamma$ is
non-normal, then $n=6$, which implies that, by GAP, $A\cong D_8\times\Bbb Z_3$ and $\Gamma$ is normal.
\eqd

Let $S$ be an inverse-closed subset of a group $G$ not containing the identity element of $G$. Recall that a permutation $\varphi$ of  $G$ is a color-preserving automorphism of $\Cay(G,S)$ if and only if we have $(xs)^\varphi\in\{x^\varphi s^{\pm 1}\}$ for each $x\in G$ and $s\in S$ \cite[p. 190]{HKMM}.

\begin{lemma}\label{rev}
 Let $R=\{a,a^{-1}\}$ and $L=\{b,b^{-1}\}$, $o(a),o(b)\geq 3$ and $R\cap L=\emptyset$. If $\Gamma$ is intransitive
then it is normal.
\end{lemma}
\proof
Suppose, towards a contradiction, that there exists $\alpha\in A$ that does not normalize $R_G$. Since $\Gamma$ is intransitive,
there is a permutation $\sigma$ of $G$ such that $(g,i)^\alpha=(g^\sigma,i)$ for all $g\in G$ and $i=1,2$. There is a natural colouring
of $\Cay(G,\{a^{\pm 1},b^{\pm 1}\})$ with two colours, where $a$-edges have one colour and $b$-edges have the other colour.
Then $\sigma$ is a
colour-preserving automorphism of $\Cay(G,\{a^{\pm 1},b^{\pm 1}\})$ because $\alpha$ is an automorphism of $\Gamma$, which means
that $(ga)^\sigma\in\{ga^{\pm 1}\}$ and $(gb)^\sigma\in\{gb^{\pm 1}\}$.

Since $\alpha$ does not normalize $R_G$ (and $G$ is $2$-generated), we know from \cite[Proposition 4.1]{HKMM} that
$G$ has a direct factor that is isomorphic to $\Bbb Z_2\times\Bbb Z_4$. So $o(a)$ and $o(b)$ are even. Therefore $o(a),o(b)\neq 3$ and so
$o(a),o(b)\geq 4$. If $o(a)=o(b)=4$ then, by GAP, $\Gamma$ is a transitive graph which is a contradiction. So we may assume that
 $o(b)>4$.

 By composing with a translation, we may assume that $\sigma$ fixes $1$. We may also assume that $\sigma$ fixes $a$ by composing
 with inversion if necessary. Then $(a^k)^\sigma=a^k$ for all $k$.

 We claim that we may assume $b^\sigma=b$. Suppose $b^\sigma\neq b$, so $b^\sigma=b^{-1}$. Then $\sigma$ is the identity
 on $\langle a\rangle$ but inverts $\langle b\rangle$, which implies that $|\langle a\rangle\cap\langle b\rangle|\leq 2$. Therefore there
 is an automorphism of $G$ agrees with $\sigma$ on $\langle a\rangle\cup\langle b\rangle$. By composing with this
 automorphism, we have $b^\sigma=b$ as desired.

 Since $\sigma$ does not normalize $R_G$, we know that $\sigma$ is not the identity permutation and so there is some minimal
 $k>0$ such that $(a^kb)^\sigma=a^kb^{-1}$. Since $a^{k-1}b$ is adjacent to $a^kb$ via an $a$-edge, we have
 $a^{k-1}ba=a^kb^{-1}$ or $a^{k-1}ba^{-1}=a^kb^{-1}$.  The first implies that $b^2=1$ which contradicts the fact
 that $o(b)\geq 4$. The second implies $a^2=b^2$. Since $|\langle a\rangle\cap\langle b\rangle|\leq 2$, we have $o(b)\leq 4$
 which contradicts the fact $o(b)>4$.
\eqd

\begin{corollary} Let $R=\{a,a^{-1}\}$ and $L=\{b,b^{-1}\}$, where $o(a)\neq o(b)$ and $R\cap L=\emptyset$.
$\Gamma$ is normal if and only if $\Gamma$ is intransitive.
\end{corollary}
\proof One direction is clear by Lemma \ref{rev}.
Let $\Gamma$ is normal and suppose, towards a contradiction, that $\Gamma$ is transitive. Then there exists $\alpha\in A$
such that $(1,1)^\alpha=(1,2)$.  Since $\Gamma$ is normal, there exists $\sigma\in\Aut(G)$
such that $\alpha=\psi_\sigma$, $R^\sigma=L$ and $L^\sigma=R$, which implies that $o(a)=o(b)$, a contradiction.
\eqd

\begin{lemma}\label{main}
Let $|R|=|L|=2$ and $R\cap L=\emptyset$. If $R=\{a,a^{-1}\}$, where $o(a)=n\geq 3$, then, perhaps after interchanging $R$ and $L$, one of the following holds:\\
$(1)$ $L=\{a^{n/2},b\}$, where $n$ is even, $b^2=1$ and $b\notin\langle a\rangle$.
In this case $G\cong\Bbb Z_n\times\Bbb Z_2$,
$\Gamma$ is
normal and intransitive.\\
$(2)$ $L=\{b,ba^{n/2}\}$, where $n$ is even, $b^2=1$ and $b\notin\langle a\rangle$. In this case $G\cong\Bbb Z_n\times\Bbb Z_2$, if $n=4$ then $\Gamma$ is transitive and non-normal and otherwise $\Gamma$ is normal and intransitive.\\
$(3)$ $L=\{b,c\}$, where $b^2=c^2=1$, $\langle a\rangle\cap\langle b,c\rangle=1$ and $b,c\notin\langle a\rangle$. In this case $G\cong\Bbb Z_n\times\Bbb Z_2^2$, if $n=4$ then $\Gamma$ is transitive and non-normal and otherwise
$\Gamma$ is normal and intransitive.\\
$(4)$ $L=\{a^k,a^{-k}\}$, for some $k\geq 2$. In this case, $G\cong\Bbb Z_n$. Furthermore, $\Gamma$
is non-normal if and only if $(n,k)=(5,2),(8,3),(10,2),(10,3),(12,5),(24,5)$. Also if $\Gamma$ is non-normal then $\Gamma$ is transitive.\\
$(5)$ $L=\{b,b^{-1}\}$, $b\notin\langle a\rangle$ and $\langle b\rangle\cap\langle a\rangle\neq 1$. In this case, $\Gamma$
is non-normal if and only if $L=\{a^3y,a^{-3}y\}$ or $L=\{a^2y,a^{-2}y\}$, where $o(a)=10$, $o(y)=2$ and 
$G=\langle a\rangle\times\langle y\rangle \cong\Bbb Z_{10}\times\Bbb Z_2$, or $L=\{ay, a^{-1}y\}$, where $o(a)=4$, $o(y)=2$ 
and $G=\langle a\rangle\times\langle y\rangle \cong\Bbb Z_4\times\Bbb Z_2$. Also if $\Gamma$ is non-normal then $\Gamma$ is transitive.\\
$(6)$ $L=\{b,b^{-1}\}$, and $\langle b\rangle\cap\langle a\rangle=1$. In this case,
$G\cong\Bbb Z_n\times\Bbb Z_l$, where $l=o(b)$  and $\Gamma$ is normal. Furthermore, $\Gamma$ is transitive if and only if $o(a)=o(b)$.
\end{lemma}
\proof It is easy to see that the only possibilities of $L$ are the cases $(1)$-$(6)$. Since $\SC(G;R,L,\{1\})\cong\SC(G;L,R,\{1\})$,
by the last three cases of Lemma \ref{2-4}, cases $(1)$, $(2)$ and $(3)$ are clear.

$(4)$ In this case, $\Gamma$ is isomorphic to the generalized Petersen graph $GP(n,k)$. Let $\Gamma$ is non-normal and
 suppose, by contrary, \[(n,k)\notin\{(5,2),(8,3),(10,2), (10,3),(12,5),(24,5)\}.\]
Then $\Gamma$ is not edge-transitive by \cite[Lemma 3 and Theorem 2]{F}.  Hence $|A|=4n$ \cite[Theorem 1 and Theorem 2]{F}. Since $\Gamma$ is non-normal, Lemma \ref{rev} and
 \cite[p. 105]{Biggs} imply that $k^2\equiv\pm 1$ (mod $n$). Hence $(k,n)=1$. Let $\sigma_1, \sigma_2,\sigma_3:G\rightarrow G$
 be the maps by the rules $(a^i)^{\sigma_1}=a^{ik}$, $(a^i)^{\sigma_2}=a^{-ik}$ and $(a^i)^{\sigma_3}=a^{-i}$. Then
 these three maps are automorphisms of $G$. Furthermore, $\psi_{\sigma_1},\psi_{\sigma_2},\varphi_{\sigma_3}\in\Aut(G;R,L)$.
 So $|\Aut(G;R,L)|\geq 4$, which implies that $A=N_A(R_G)$ i.e $\Gamma$ is normal, a contradiction.

 Conversely, suppose that \[(n,k)\in\{(5,2),(8,3),(10,2),(10,3),(12,5),(24,5)\}.\] Then $\Gamma$ is arc-transitive by \cite[p. 105]{Biggs} and
 so it is non-normal by Lemma \ref{arc}. If
 $\Gamma$ is non-normal then it is transitive by Lemma \ref{rev}.

 $(5)$ Let $L=\{b,b^{-1}\}$, $b\notin\langle a\rangle$, $a\notin\langle b\rangle$ and $\langle b\rangle\cap\langle a\rangle\neq 1$. If $\Gamma$ is intransitive, then by Lemma \ref{rev}, $\Gamma$ is normal. Hence, we may assume that $\Gamma$ is transitive. If $\Gamma$ is arc-transitive then by \cite[Proposition 5.1(2)]{Zhou-Feng}, it is  the unique arc-transitive cubic graph of order 40, denoted by F040A in the Foster Census,
 $G=\langle x\rangle\times\langle y\rangle\cong\Bbb Z_{10}\times\Bbb Z_{2}$, and we may assume that $a=x$ and $b\in\{x^3y,x^2y\}$.  Then $\Gamma$ is non-normal by Lemma \ref{arc}. So we may now assume that $\Gamma$ is not arc-transitive. Then by \cite[Theorem 1.1]{Zhou-Feng}, up to isomorphism, one of the following happens:\\
$(i)$
 $G=\langle x\rangle\times\langle y\rangle\cong\Bbb Z_{mk}\times\Bbb Z_m$, $k\geq 3,m\geq 1$, where $(m,k,t)=(1,10,2)$ or 
 $(t,mk)=1$ and
 $t^2\equiv -1$ (mod $k$) and we may assume that $a=x$, $b=x^ty$.  Clearly $m=1$ is impossible, because $b\notin\langle a\rangle$. Also $\Aut(\Gamma)\cong R_G\rtimes\Bbb Z_4$ \cite[Theorem 5.5(3)]{Zhou-Feng}, which implies that $\Gamma$ is normal.\\
 $(ii)$
 $G=\langle x\rangle\times\langle y\rangle\cong\Bbb Z_{mk}\times\Bbb Z_m$, $km\geq 3$ and $m\geq 1$, where
 $(t,mk)=1$, $t^2\equiv 1$ (mod $k$) and $\Gamma$ is a Cayley graph over $G\rtimes\langle z\rangle$ for some
 involution $z$,  \cite[Theorem 5.2(5)]{Zhou-Feng}.
 Furthermore, we may assume that $a=x$, $b=x^ty$. Clearly $m=1$ is impossible because $b\notin\langle a\rangle$. 
 
Since $\Gamma$ is connected and transitive but not edge-transitive, every automorphism of $\Gamma$ maps $G$-orbits to $G$-orbits.
If $(m,k)=(2,2)$, then $L=\{ay,a^{-1}y\}$,
$G\cong\Bbb Z_4\times\Bbb Z_2$ and $\Gamma$
is non-normal over $G$ by GAP. Hence, we may assume that $(m,k)\neq (2,2)$. Then we claim that $\Gamma$ is normal. Suppose,
towards a contradiction, that $\Gamma$ is not normal. So, there exists a color-preserving automorphism $\sigma$ of $\Gamma_0=\Cay(G,\{a,a^{-1},b,b^{-1}\})$ which fixes $1$ but is not a group automorphism of $G$ (see \cite[Remark 2.1]{HKMM}). Since the map $x\mapsto x^{-1}$ is an automorphism of $G$, we may assume that $a^\sigma=a$. 

We may assume that $\sigma$ is not the identity. Then there is some
$i$ such that $(a^ib)^\sigma=a^ib^{-1}$. By composing with a translation, we may assume that $i=0$ and $b^\sigma=b^{-1}$. Then we have $(b^m)^\sigma=b^{-m}$. But $b^m\in\langle a\rangle$ and $\sigma$ is the identity on $\langle a\rangle$. So $b^m$
must have order two which means that
$b^m=a^{km/2}$. So $tm\equiv km/2$ (mod $km$) which means that $t\equiv k/2$ (mod $k$). Since $t^2\equiv 1$ (mod $k$) this implies that $k=2$. So, by \cite[Corollary 4.2]{HKMM},  $|G|=2m^2$ is divisible by $8$. Thus $m$ is even and $m\geq 4$ because $(m,k)\neq (2,2)$. Since $\langle a\rangle\cap\langle b\rangle\neq 1$ and $(b^m)^{-1}=b^m$, the 
map $\varphi:a^ib^j\mapsto a^ib^{-j}$ is a well-defined automorphism of $G$ that is also
an automorphism of $\Gamma_0$. Furthermore, $\psi=\sigma\varphi$
is a color-preserving automorphism 
of $\Gamma_0$ which fixes all powers of $a$ (including $1$) and $b$, but is not a group automorphism of $G$.

Since $mk,m\geq 4$, it is easy to see that for all $g\in G$, $g$ and $gab$ are the only common neighbours of $ga$ and $gb$ in $\Gamma_0$. Putting $g=1$, we get
$(ab)^\psi=ab$. Now putting $g=a$ we get $(a^2b)^\psi=a^2b$. By continuing this procedure we get $(a^ib)^\psi=a^ib$ for all $i$. Since $m\geq 4$, we have $b^2\neq 1$.  So, for all $i$ we have $(a^ib^2)^\psi=a^ib^2$. This implies that $(a^ib^3)^\psi=a^ib^3$ for all $i$. By continuing this procedure,
we get $(a^ib^j)^\psi=a^ib^j$ for all $i,j$. This means that $\psi$ is the trivial automorphism of $G$. Hence $\sigma$ is an automorphism of group $G$, a contradiction.

(6) Let $L=\{b,b^{-1}\}$, where 
 $\langle a\rangle\cap\langle b\rangle=1$. We claim that $\Gamma$ is normal. If $\Gamma$ is intransitive, then by Lemma \ref{rev}, $\Gamma$ is normal. Hence, we may assume that $\Gamma$ is transitive.
 Then, by \cite[Theorem 1.1]{Zhou-Feng},
 $G=\langle x\rangle\times\langle y\rangle$, $o(x)=mk$, $o(y)=m$, for some $m,k\geq 1$, where $mk\geq 3$.
 Furthermore 
 $a=x$
 and $b=a^ty$ for some integer $t$ with $(t,mk)=1$ and
 $t^2\equiv 1$ (mod $k$), or $(m,k,t)=(1,10,2)$, or
 $(t,mk)=1$ and $t^2\equiv -1$ (mod $k$). Clearly $(m,k,t)=(1,10,2)$ is impossible, because $\langle a\rangle\cap\langle b\rangle=1$.
So we have $b^m=a^{tm}\in\langle a\rangle\cap\langle b\rangle=1$. Thus $o(b)$ divides $m$, and $k$ divides
$t$. The latter implies that $k=t=1$. Thus $b=ay$ and
$o(b)=o(a)=m$. 

Since $\Gamma$ is connected, \cite[Proposition 5.1]{Zhou-Feng} implies that
$\Gamma$ is not edge-transitive. So every automorphism of $\Gamma$ maps $G$-orbits to $G$-orbits. Suppose, towards a contradiction, that
$\Gamma$ is not normal. Similar to the previous case,
there exists a color-preserving automorphism $\sigma$ of $\Gamma_0=\Cay(G,\{a,a^{-1},b,b^{-1}\})$ which fixes $1$ but is not a group autmorphism of $G$ and we may assume that $a^\sigma=a$ and $b^\sigma=b^{-1}$. Then \cite[Theorem 1.3(ii)]{HKMM} implies that $8$ divides $|G|=m^2$.
So $4$ divides $m$. Since $\langle a\rangle \cap\langle b\rangle=1$, $\varphi: a^ib^j\mapsto a^ib^{-j}$ is a well-defined automorphism of $G$ that is also an automorphism of $\Gamma_0$. Again, by the same argument in the last paragraph of the proof of previous case, we get $\sigma$ is an automorphism of $G$ which is a contradiction. So we have proved that
$\Gamma$ is normal.

As we saw above, if $\Gamma$ is transitive, then $o(a)=o(b)$. Conversely, suppose that $o(a)=o(b)$. Then $\sigma:a^ib^j\mapsto a^jb^i$ is a group automorphism of $G$ and $\langle R_G,\psi_\sigma\rangle$, where $\psi$ is defined by the rule $(g,1)^\psi=(g^\sigma,2),(g,2)^\psi=(g^\sigma,1)$ for all
$g\in G$, is a transitive subgroup of $\Aut(\Gamma)$.
This completes the proof.\eqd

\textbf{Proof of Theorem \ref{asli}} It is a direct consequence of Lemmas \ref{3.3}, \ref{3.4}, \ref{3.6}-\ref{2-4} and \ref{main}.\eqd

\textbf{Acknowledgements}
The authors gratefully appreciate anonymous referee for constructive comments and recommendations which
definitely helped to improve the readability and quality of the paper.



\end{document}